\documentclass[conference]{IEEEtran}
\IEEEoverridecommandlockouts
\usepackage{cite}
\usepackage{amsmath,amssymb,amsfonts}
\usepackage{bm} 
\usepackage[shortlabels]{enumitem}  
\usepackage{algorithmic}
\usepackage{graphicx}
\usepackage{textcomp}
\usepackage{xcolor}
\usepackage[linesnumbered,ruled,vlined]{algorithm2e} 
\usepackage{mathrsfs} 
\def\BibTeX{{\rm B\kern-.05em{\sc i\kern-.025em b}\kern-.08em
    T\kern-.1667em\lower.7ex\hbox{E}\kern-.125emX}}

\begin{document}
\title{An Efficient Algorithm to Generate all Labeled Triangle-free Graphs with a given Graphical Degree Sequence*\\
\thanks{ 
This research has been supported by a research seed grant of the College of Engineering and
Computing at Georgia Southern University. The computational experiments have been supported by the Talon
 Cluster at Georgia Southern University.
}
}

\author{\IEEEauthorblockN{Kai Wang}
\IEEEauthorblockA{\textit{Department of Computer Science} \\
\textit{Georgia Southern University}\\
Statesboro, GA, USA \\
kwang@georgiasouthern.edu}
}

\maketitle

\begin{abstract}
We extend our previous algorithm that generates all labeled graphs with a given graphical degree sequence
to generate all labeled triangle-free graphs with a given graphical degree sequence. The algorithm uses
various pruning techniques to avoid having to first generate all labeled realizations of the input sequence
and then testing whether each labeled realization is triangle-free. It can be further extended to
generate all labeled bipartite graphs with a given graphical degree sequence by adding a simple test whether
each generated triangle-free realization is a bipartite graph. All output graphs are generated in the
lexicographical ordering as in the original algorithm. The algorithms can also be easily parallelized.
\end{abstract}


\section{Introduction}
A finite simple graph $G=(V,E)$ is called a triangle-free graph if it contains no three vertices forming
a clique $K_3$. In general, it is called $K_r$-free if the vertex set $V$ contains no $r$-subset forming a
clique $K_r$. Given a graph $G$ of order $n$, it takes at most $O(n^3)$ time to decide
whether it is triangle-free. However, it is NP-complete to decide if it contains a clique $K_r$, where $r$
is an input parameter together with $G$ satisfying $3\le r\le n$ \cite{Karp1972}.

A related decision problem $P(\bm{d},r)$ is as follows. Given a graphical degree sequence $\bm{d}=\{d_1\ge d_2\ge
\cdots\ge d_n\}$ with $d_n\ge 1$, and a parameter $r$ such that $3\le r\le n$, is there a simple graph realization of
$\bm{d}$ that is $K_r$-free? In the special case $r=3$, the problem is to decide if $\bm{d}$ has a triangle-free
realization. Although it is conjectured this general decision problem $P(\bm{d},r)$ is intractable \cite{Bauer2008}, we
know of no rigorous proof that this problem is NP-complete. In fact, we conjecture that this problem is not
in NP because the formal language $\overline{L}=\{\langle G,r\rangle|\mbox{ The graph } G \mbox{ does not
have a clique of size } r\}$ (which is the complement of the NP-complete formal language
$L=\{\langle G,r\rangle|\mbox{ The graph } G\mbox{ has a clique of size } r\}$) 
is coNP-complete (using the fact that the complement of an NP-complete language is coNP-complete)
and every coNP-complete language is conjectured to be outside of NP \cite{Goldreich2008}.
Nevertheless, the special case problem $P(\bm{d},3)$ is in NP since it takes at most
$O(n^3)$ time to decide whether a given graph is triangle-free so that a triangle-free realization of $\bm{d}$
(if it exists) would be a certificate for membership of any input $\bm{d}$ to be in the formal language $P(\bm{d},3)
=\{\langle \bm{d}\rangle|\mbox{ The graphical sequence } \bm{d}\mbox{ has a $K_3$-free realization.}\}$.
However, to obtain a ``yes'' answer to this special-case decision problem deterministically such a certificate may
not be necessary.

One may also consider the following combinatorial generation problems related to the decision problem $P(\bm{d},3)$:
\begin{enumerate}[(A)]
	\item Generate some labeled triangle-free realization of $\bm{d}$.
	\item Generate all labeled triangle-free realizations $G=(V,E)$ of $\bm{d}$. Denote the set of all labeled triangle-free
	realizations of $\bm{d}$ as $\mathcal{S}(\bm{d})$.
	\item Generate a uniformly random labeled triangle-free realization of $\bm{d}$.
\end{enumerate}
These problems appear to be more sophisticated than the corresponding combinatorial generation problems in
\cite{Wang2024} where we seek realizations of $\bm{d}$ without imposing any additional structural restrictions.
In this paper we will mainly deal with the problem (B).

A special subclass of triangle-free graphs is bipartite graphs. There are triangle-free graphs of arbitrarily high chromatic
number \cite{Mycielski1955} and it is NP-complete to decide if a given triangle-free graph $G$ of order $n$ has chromatic
number $\le r$, where $r$ is an input parameter together with the graph $G$ satisfying $3\le r\le n$ \cite{Maffray1996}.
Since it is well-known that it is easy to decide if a (triangle-free or not) graph is bipartite, we can also consider variations
of the above three combinatorial generation problems where labeled bipartite realizations of $\bm{d}$ are to be generated.

The rest of the paper is organized as follows. In Section \ref{sec:review_alg} we review our basic generation
algorithm to systematically generate all labeled realizations of a given graphical sequence in \cite{Wang2024}.
In Section \ref{sec:our_alg_tf} and \ref{sec:our_alg_bp} we present our simple algorithms to generate all labeled
triangle-free and bipartite realizations of a given graphical sequence, respectively. In Section \ref{sec:analysis}
we provide a rough complexity analysis of our algorithms. In Section \ref{sec:experiments} we show some
experimental results. In Section \ref{sec:conclusion} we conclude with some future research directions.

\section{Review of the Basic Generation Algorithm}
\label{sec:review_alg}
The overall framework of our exhaustive generation algorithm for all labeled triangle-free realizations of an input $\bm{d}$
is similar to that presented in \cite{Wang2024} and it will be
a backtracking algorithm where each node of the search tree $\mathcal{T}(\bm{d})$ traversed and built by the
algorithm represents a partial or finished labeled triangle-free realization of the input sequence $\bm{d}$. Each finished
realization $G$ has the same vertex set $V=\{v_1,v_2,\cdots,v_n\}$ satisfying (i) $d_G(v_i)=d_i$ for $i=1,\cdots,n$;
and (ii) $G$ has no triangles. The total ordering of all these labeled triangle-free realizations is the same derived
lexicographical ordering of labeled realizations based on the lexicographical ordering of the vertex labels
$v_1<v_2<\cdots<v_n$ (see \cite{Wang2024} for the definition of this total ordering).
Our algorithm will generate all output triangle-free graphs based on this ordering.
We will still use the terminology ``node'' to mean a point in the search tree $\mathcal{T}(\bm{d})$ and use the terminology
``vertex'' to mean a point in a labeled triangle-free realization $G$ of $\bm{d}$ to be generated by the algorithm.
At any moment during the execution of the algorithm the \textit{residual degree} of a vertex $v$ in $G$ is a dynamic quantity
representing the number of neighbors yet to be assigned for $v$ to fulfill its degree requirement specified in
$\bm{d}$. A vertex of $G$ is said to be \textit{saturated} or \textit{unsaturated} based on whether its residual degree is
0 or not. For comparison and descriptive purposes, we use $\mathcal{T}_0(\bm{d})$ to denote the search tree traversed
and built by the general exhaustive generation algorithm for all labeled realizations in \cite{Wang2024} while we use
$\mathcal{T}(\bm{d})$ to denote the search tree traversed and built by the
algorithm for triangle-free realizations in this paper.

\section{Our Generation Algorithm for Triangle-free Realizations}
\label{sec:our_alg_tf}
We assume the input $\bm{d}$ is a graphical degree sequence. (There are several equivalent and efficient criteria
to decide whether this is the case \cite{Sierksma1991}, among which the Erd{\H{o}}s-Gallai criterion \cite{ErdosCallai1960}
is probably the most well-known. The Havel-Hakimi criterion \cite{Havel1955,Hakimi1962} is also well-known for accomplishing job.)
However, it is not assumed to have a triangle-free
realization. Thus, before our algorithm starts to build the search tree $\mathcal{T}(\bm{d})$ to construct all triangle-free
realizations, it runs three simple checks (we call these checks stage 1 operations) to rule out the easy cases
where $\bm{d}$ has no triangle-free realizations.

By Mantel's theorem \cite{Mantel1907}, the maximum number of edges in a triangle-free simple graph of order $n$ is
$\lfloor \frac{n^2}{4}\rfloor$. Thus, if the sum of degrees $S=\sum_{i=1}^{n}d_i$ of the input $\bm{d}$ is such that
$S>\lfloor \frac{n^2}{2}\rfloor$ (recall that the degree sum of a simple graph is twice the number of its edges), then our
algorithm can immediately report that there are no triangle-free realizations and stop.

The next rule to show that the input $\bm{d}$ has no triangle-free realizations uses the concept of
\textit{residue} of an arbitrary graphical degree sequence introduced in \cite{Favaron1991}. The residue $R(\bm{d})$ of
a graphical degree sequence $\bm{d}=\{d_1\ge d_2 \ge \cdots \ge d_n\}$ is the number of
zeros obtained in the last step by the iterative Havel-Hakimi reduction procedure \cite{Havel1955,Hakimi1962}.
As shown in \cite{Favaron1991}, the residue $R(\bm{d})$ of a graphical degree sequence $\bm{d}$ is a lower bound
on the independence number of every realization of $\bm{d}$. If we denote by $\overline{\bm{d}}$ the complementary
graphical degree sequence of $\bm{d}$ (i.e. $\overline{\bm{d}}=\{n-1-d_n\ge n-1-d_{n-1}\ge \cdots \ge n-1-d_1\}$),
which is the degree
sequence of the complementary graph of any realization of $\bm{d}$, then our algorithm can immediately report that there
are no triangle-free realizations and stop when it finds that $R(\overline{\bm{d}})\ge 3$. This is because when this
is the case every realization of $\bm{d}$ will have a clique number $\ge 3$ (i.e. contain a triangle). We note that $R(\bm{d})$
can be computed easily in $O(n^2)$ time for any input $\bm{d}$ of length $n$. For our algorithm most interesting inputs $\bm{d}$
satisfy $R(\overline{\bm{d}})=2$ ($R(\overline{\bm{d}})=1$ if and only if $\bm{d}$ is the degree sequence of $n$ zeros,
which cannot happen when we assume $d_n\ge 1$ for our input $\bm{d}$). Also, $R(\bm{d})$ is not necessarily a tight
lower bound on the minimum independence number among all realizations of $\bm{d}$. In fact, the difference between the
two quantities can become arbitrarily large \cite{Jelen1996}. Thus, when our algoirhtm finds out that
$R(\overline{\bm{d}})=2$ it will proceed to stage 2 (see below) with a possibility of reporting no triangle-free realizations in the end.

The final simple rule our algorithm adopts to show that the input $\bm{d}$ has no triangle-free realizations
uses the concept of \textit{Murphy's bound} introduced in Murphy \cite{Murphy1991},
denoted $\beta(\bm{d})$ here, which is also a lower bound on the independence number
of any realization of $\bm{d}$. Therefore, our algorithm can also immediately report that there
are no triangle-free realizations and stop when it finds that $\beta(\overline{\bm{d}})\ge 3$. Murphy's bound
$\beta(\bm{d})$ can be computed easily in $O(n)$ time for any given input $\bm{d}$ of length $n$. Similar to the situation
of residue, most interesting inputs $\bm{d}$ satisfy $\beta(\overline{\bm{d}})=2$ ($\beta(\overline{\bm{d}})=1$
if and only if $\bm{d}$ is the degree sequence of $n$ zeros). Also, $\beta(\bm{d})$ is not necessarily a tight lower bound on the
minimum independence number among all realizations of $\bm{d}$, which forces our algorithm to proceed into
stage 2 when it finds out that $\beta(\overline{\bm{d}})=2$ (after finding out $R(\overline{\bm{d}})=2$). Note that there exists
$\bm{d}$ such that $\beta(\overline{\bm{d}})\ge 3$ and $R(\overline{\bm{d}})=2$ and there also exists $\bm{d}$
such that $R(\overline{\bm{d}})\ge 3$ and $\beta(\overline{\bm{d}})=2$.

The above three simple rules are actually also applied in deciding whether a given input $\bm{d}$ has
a bipartite realization in \cite{Wang2021}. If an input $\bm{d}$ passes the tests of the above three
rules in stage 1, then our algorithm will enter stage 2 to build the search tree $\mathcal{T}(\bm{d})$ as in
\cite{Wang2024} to construct all labeled triangle-free realizations. We emphasize it is possible that an input $\bm{d}$ can
pass the three rule checks in stage 1 and still have no triangle-free realizations. This means that the algorithm is possible
to do a lot of work in building and traversing the search tree $\mathcal{T}(\bm{d})$ and then
report that no triangle-free realizations are found at the end of stage 2. 
We note that the complexity of the decision version of the problem $P(\bm{d},3)$ is unknown and we suspect that
there is a chance that the decision problem is not polynomial-time solvable. Therefore, our current generation
algorithm is not always able to avoid stage 2 in case the input $\bm{d}$ has no triangle-free realizations.

In the second stage the algorithm starts at the root node of the search tree $\mathcal{T}(\bm{d})$, which represents
a labeled graph with a set $V$ of isolated vertices $\{v_1,v_2,\cdots,v_n\}$. In general, for each $1\le k\le n$, at depth
$k$ of the search tree $\mathcal{T}(\bm{d})$ are those search tree nodes representing the partial labeled graphs $G$
being constructed in which the vertices $v_1,v_2,\cdots,v_k$ are already saturated (i.e. they have
$d_1,d_2,\cdots,d_k$ neighbors already assigned, respectively) with no triangles in this partial realization $G$ so far
and the remaining vertices $v_{k+1},v_{k+2},\cdots,v_n$ can potentially have additional neighbors properly added to
produce a finished labeled triangle-free realization of $\bm{d}$ if there is a labeled triangle-free realization of $\bm{d}$
containing as a subgroup the current partial realization $G$ with the current neighbor assignments for $v_1,v_2,\cdots,v_k$.


As in \cite{Wang2024}, our algorithm will traverse the search tree nodes of $\mathcal{T}(\bm{d})$ in a depth-first
manner and only one copy of the (partial or finished) labeled triangle-free realization $G$ of $\bm{d}$ will be
maintained throughout the traversal of the entire search tree and the
algorithm will dynamically update this copy of labeled graph $G$ when traversing from one search tree node to the next.
When traversing a search tree node at depth $k$, the algorithm behaves differently based on whether the vertex
$v_k$ in the labeled graph $G$ is already saturated or not.

If $v_k$ is not saturated, then the algorithm will build a candidate set of neighbors $C_k$ for $v_k$, where $C_k$ will be
a pruned subset of the set of those vertices in $U_k=\{v_{k+1},\cdots,v_n\}$ that are not already
saturated in the copy of our labeled triangle-free graph $G$ being constructed. We do not need to consider any vertex
in $L_k=\{v_1,v_2,\cdots,v_{k-1}\}$ for inclusion into $C_k$ because at depth $k$ of the search tree all vertices in
$L_k$ must have already been saturated. The following pruning will be applied in constructing the set $C_k$:
\begin{enumerate}[(1)]
	\item if an unsaturated vertex $y\in U_k$ is adjacent to any neighbor $x\in L_k$ of $v_k$, then
	it will not be included in $C_k$. This is because we want to avoid a triangle in $G$ consisting of $\{v_k,x,y\}$.
\end{enumerate}
If the above pruning is not applied, we would always have $|C_k|\ge r_k$, where $r_k$ is the residual degree of $v_k$
before we add any neighbors for $v_k$ in $G$ at the current search tree node, because the input $\bm{d}$ is assumed
to be a graphical degree sequence. When the pruning is applied, it is possible that $|C_k|< r_k$ so that $v_k$ cannot be
saturated in a triangle-free manner and this part of
the search tree with root at the current search tree node that would appear in $\mathcal{T}_0(\bm{d})$ is pruned.
If the pruned candidate set of neighbors $C_k$ for $v_k$ does contain at least $r_k$ elements, our algorithm will then
use the subset enumeration algorithm to select the lexicographically smallest set of $r_k$ neighbors from $C_k$
such that adding this set of $r_k$
neighbors for $v_k$ in $G$ will not make it impossible to finish a labeled triangle-free realization $G$ of $\bm{d}$
by properly adding additional neighbors for $\{v_{k+1},\cdots,v_n\}$. 
As in \cite{Wang2024}, the algorithm will also need to verify that the sequence $\bm{r}$ of residual degrees of the
$n$ vertices resulting from adding a possible $r_k$-subset of neighbors from $C_k$ for $v_k$ is still graphical using
any of several equivalent criteria to decide whether a given integer sequence is a graphical sequence \cite{Sierksma1991}.
However, before (or after) it performs the graphicness check for each sequence of residual degrees resulting from
a selection of an $r_k$-subset of neighbors from $C_k$ for $v_k$, the algorithm performs three additional checks
shown below to determine whether the current
selection of $r_k$-subset of neighbors from $C_k$ for $v_k$ would make it impossible to lead to a finished
labeled triangle-free realization $G$ of $\bm{d}$.
\begin{enumerate}[label=(\roman*)]
	\item Let $\sigma_1$ be the sum of residual degrees of those neighbors in the selected $r_k$-subset  from $C_k$
	for $v_k$ after their projected adjacency to $v_k$. Let $\sigma_2$ be the sum  of residual degrees of those vertices
	with positive residual degrees but are not in the currently selected $r_k$-subset of $C_k$. Note that these vertices for
	computing $\sigma_2$ may or may not be in $C_k$ since certain vertices with positive residual degrees may not have
	been added into $C_k$ in order to avoid a triangle in $G$ (see our previous pruning when constructing $C_k$). If
	$\sigma_1>\sigma_2$, then the currently selected $r_k$-subset
	of $C_k$ is infeasible to saturate $v_k$ in the sense that it would make it impossible to finish a labeled triangle-free
	realization $G$ of $\bm{d}$. This is because these $r_k$ vertices from $C_k$ would all be made neighbors of $v_k$
	and we cannot make any two of them adjacent in any decendent search tree nodes of the current search tree node
	in order to avoid a triangle in $G$. This means further neighbors added for each of these $r_k$ vertices must come
	from those vertices still with positive residual degrees but are outside of the current $r_k$-subset of neighbors from $C_k$
	for $v_k$. Thus, if $\sigma_2<\sigma_1$, then there is no way to saturate all vertices in the current $r_k$-subset of
	neighbors for $v_k$ in a triangle-free manner while also keeping the currently selected neighbors for $v_k$ in the
	decendent search tree nodes.
	\item when computing $\sigma_1$, we keep track of the maximum residual degree $m$ of those $r_k$ vertices selected
	from $C_k$ as neighbors for $v_k$ after they had all been made adjacent to $v_k$ on the current search tree node. When
	computing $\sigma_2$, we keep track of the number $c$ of vertices that are still with positive residual degrees but are not
	in the currently selected $r_k$-subset of $C_k$ for $v_k$ (recall $\sigma_2$ is the sum of residual degrees of these
	$c$ vertices: see above). If $m>c$, then the currently selected $r_k$-subset of $C_k$ is also infeasible to saturate $v_k$.
	As indicated above, any future neighbors for each of these $r_k$ selected neighbors from $C_k$ for $v_k$ in any
	decendent search tree node of the current search tree node must come from the set of those vertices outside of
	the current $r_k$-subset of $C_k$ but with positive residual degrees in order to avoid a triangle in $G$. If $c<m$,
	then the vertex in the $r_k$-subset of $C_k$ with the maximum residual degree $m$ (recall that $m$ is computed
	after the vertex had been made adjacent to $v_k$) would not have enough number of neighbors for it in any decendent
	search tree node.
	\item if $m=c$ while $\sigma_1\ne \sigma_2$, then the currently selected $r_k$-subset of $C_k$ is still infeasible to
	saturate $v_k$. This is because all those $c$ vertices with positive residual degrees outside of the selected
	$r_k$-subset of $C_k$, after they had \textit{all} been made adjacent to the vertex with maximum residual
	degree $m$ in the selected $r_k$-subset of $C_k$, cannot have additional edges added between any two
	of them (again, to avoid a triangle in $G$). Therefore, any future neighbors for each of those $c$ vertices in the
	decendent search tree nodes must come from the currently selected $r_k$-subset of $C_k$.
	Thus, it is required that $\sigma_1-m=\sigma_2-c$ for feasibility (i.e. $\sigma_1=\sigma_2$).
\end{enumerate}
The above three feasibility checks (i) $\sigma_1\le\sigma_2$ (ii) $m\le c$ and (iii) $m=c\Rightarrow \sigma_1=\sigma_2$
for each enumerated $r_k$-subset of $C_k$ can be performed either before or after checking the graphicness of the sequence
of residual degrees $\bm{r}$ resulting from adding the $r_k$-subset of neighbors from $C_k$ for $v_k$. If the three
feasibility checks did not all pass, then the graphicness test of $\bm{r}$ can be skipped. If the graphicness test of $\bm{r}$
did not pass, then the three feasibility checks can all be skipped. Once the algorithm finds the lexicographically smallest
subset of $r_k$ neighbors from $C_k$ for $v_k$ that passes all the three feasibility checks and the graphicness test of $\bm{r}$,
it just updates the copy of the labeled graph $G$ being constructed accordingly and then traverses downward to
the next depth $k+1$ of the search tree $\mathcal{T}(\bm{d})$. As mentioned in \cite{Wang2024}, we would need a
fundamental combinatorial algorithm to
generate all $K$-subsets of a set of $M$ elements in lexicographical order for the purpose of finding a potentially
feasible subset of $r_k$ neighbors from $C_k$ for $v_k$ (with $K=r_k,M=|C_k|$ here).
We note that passing all these tests and finding a potentially feasible subset of $r_k$ neighbors from $C_k$ for $v_k$
does not guarantee that there are labeled triangle-free realizations of $\bm{d}$ in the part of the search tree $\mathcal{T}(\bm{d})$
rooted at the current search tree node.

We have finished the description of the behavior of our generation algorithm for the situaiton when $v_k$ is not saturated.
If $v_k$ is already saturated, then the algorithm will behave similarly to that in \cite{Wang2024} under the same situation.
Let the set of neighbors of the saturated vertex $v_k$ in the copy of the labeled graph $G$ be $N(v_k)=L(v_k)\cup R(v_k)$ (where
$L(v_k)$ is the set of neighbors of $v_k$ in $G$ from $L_k=\{v_1,\cdots,v_{k-1}\}$ and $R(v_k)$ is the set of neighbors
of $v_k$ in $G$ from $U_k=\{v_{k+1},\cdots,v_n\}$), then the algorithm will directly proceed to traverse downward to
depth $k+1$ in the search tree $\mathcal{T}(\bm{d})$ if it is visiting the current search tree node for the first time. If it is
not visiting the current search tree node for the first time (i.e. it is visiting the node from a downward search tree node
at depth $k+1$ during a backtracking), then it will try to find a new set $N'(v_k)=L(v_k)\cup R'(v_k)$ of neighbors for
$v_k$ by finding the next valid lexicographically larger set of neighbors $R'(v_k)$ for $v_k$ (i.e. those neighbors that
are lexicographically larger than $v_k$) with 
$|R'(v_k)|=|R(v_k)|$ that can potentially lead to a finished labeled triangle-free realization $G$ of $\bm{d}$. In order to
find $R'(v_k)$ for $v_k$, we do not need to enumerate all lexicographically larger $|R(v_k)|$-subsets of $U_k=
\{v_{k+1},\cdots,v_n\}$ than current $R(v_k)$. Instead, we can first form a forbidden subset of neighbors $F_k$ for $v_k$
consisting of all those vertices $y$ in $U_k$ that are adjacent to any neighbor $x$ of $v_k$ in $L_k=\{v_1,\cdots,v_{k-1}\}$.
That is, $F_k=\{y\in U_k|yx\in E(G)\mbox{ for some }x\in L(v_k)\}$ and this is similar to our previous pruning technique
to avoid triangles in $G$ when constructing $C_k$ in considering the situation that $v_k$ is not saturated.
After $F_k$ is formed, our algorithm then enumerates all lexicographically larger $|R(v_k)|$-subsets of $U_k-F_k$
than current $R(v_k)$. This potentially helps reduce the number of enumerated subsets in finding $R'(v_k)$ and
prune many search tree nodes that would appear in $\mathcal{T}_0(\bm{d})$. For each
enumerated $|R(v_k)|$-subset of $U_k-F_k$, we would still need to perform a graphicness test of the
sequence of residual degrees $\bm{r'}$ constructed in the same way as in \cite{Wang2024}. If
the next valid lexicographically larger set of neighbors $R'(v_k)$ for $v_k$ that passes the graphicness test
of $\bm{r'}$ can be found, then the copy of the labeled graph $G$ being constructed is actually
updated accordingly and the algorithm proceeds downward to depth $k+1$ in the search tree.
If no such set $R'(v_k)$ is found for $v_k$, then the algorithm backtracks to the parent node of the current
search tree node without making any changes to the copy of the labeled graph $G$ being constructed.
Again, here we would also need the fundamental combinatorial algorithm to generate $K$-subsets of a set of
$M$ elements (with $K=d_k-|L(v_k)|=|R(v_k)|,M=|U_k-F_k|$) that are lexicographically larger than $R(v_k)$ for
the purpose of finding $R'(v_k)$. Also notice that here we do not need any feasibility checks like the above
feasibility checks (i)-(iii) in considering the situation when $v_k$ is not saturated because of the introduction
of the forbidden set $F_k$ before trying to find $R'(v_k)$.

Interestingly, we can treat the search tree $\mathcal{T}(\bm{d})$ built and traversed by our generation algorithm
as a subtree of the search tree $\mathcal{T}_0(\bm{d})$ built when trying to construct all labeled realizations of
the input graphical sequence $\bm{d}$
as in \cite{Wang2024}. If the input $\bm{d}$ has no triangle-free realizations, then all leaf nodes of
$\mathcal{T}_0(\bm{d})$ are pruned and become absent in $\mathcal{T}(\bm{d})$, and this conclusion will be
reported when the algorithm finishes traversing $\mathcal{T}(\bm{d})$ and terminates at the end of stage 2.
For most input sequences $\bm{d}$ that have triangle-free realizations
our algorithm will retain exactly those leaf nodes in $\mathcal{T}_0(\bm{d})$ representing the triangle-free realizations
while having the other leaf nodes (representing labeled realizations with clique number $\ge 3$) pruned somewhere in
the internal nodes of the search tree. In fact, the largest clique number $\Omega(\bm{d})$ of any realization of the input
$\bm{d}$ can be efficiently computed \cite{Yin2012}, in contrast to the intractabilithy of computing the minimum
clique number $\omega(\bm{d})$ of any realization of the input $\bm{d}$ as conjectured in \cite{Bauer2008}.
Due to pruning, the overall runtime of traversing and building $\mathcal{T}(\bm{d})$ for many inputs $\bm{d}$ can
be much shorter than traversing and building the entire $\mathcal{T}_0(\bm{d})$, as demonstrated later in Section
\ref{sec:experiments}. For certain rare input sequences $\bm{d}$, all their realizations are triangle-free. In this
case $\mathcal{T}(\bm{d})$ and $\mathcal{T}_0(\bm{d})$ will be the same (i.e. no actual prunings happen).
Clearly, if $\Omega(\bm{d})=2$, then $\mathcal{T}(\bm{d})$ and $\mathcal{T}_0(\bm{d})$ will be identical.
If $\Omega(\bm{d})>2$, then $\mathcal{T}(\bm{d})$ will be a proper subtree of $\mathcal{T}_0(\bm{d})$.
If $\omega(\bm{d})=2$, then $\mathcal{T}(\bm{d})$ will contain some leaf nodes of $\mathcal{T}_0(\bm{d})$.
If $\omega(\bm{d})>2$, then $\mathcal{T}(\bm{d})$ will contain no leaf nodes of $\mathcal{T}_0(\bm{d})$.

Finally we note that the stage 2 of building and traversing the search tree $\mathcal{T}(\bm{d})$ can be parallelized
in a similar way as that of parallelizing the algorithm of building and traversing the search tree $\mathcal{T}_0(\bm{d})$
in \cite{Wang2024}, though the load balancing problem of the parallel algorithm would still be challenging.

\section{Extension Generation Algorithm for Bipartite Realizations}
\label{sec:our_alg_bp}
A closely related combinatorial generation problem is to generate all bipartite realizations of an input graphical sequence
$\bm{d}$. Since any simple graph of order $n$ can be easily tested whether it is bipartite (or, equivalently, 2-colorable)
in $O(n^2)$ time, we can further extend our previous generation algorithm in Section \ref{sec:our_alg_tf} by simply adding a test
whether each generated triangle-free realization $G$ is bipartite. In this way all labeled bipartite realizations of the input
sequence $\bm{d}$ can also be generated in lexicographical order. Also, before the extension algorithm enters stage 2
to build and traverse the search tree $\mathcal{T}(\bm{d})$, it can actually employ several additional simple and
efficient checks in \cite{Wang2021} (a total of seven rules are listed there) as stage 1 operations to rule out the cases
where the input $\bm{d}$ does not have any bipartite realizations. In fact the whole decision algorithm there for
potential bipartiteness of the input $\bm{d}$ can be applied here for our extension generation algorithm during stage 1
before entering stage 2, with stage 2 skipped if stage 1 decides that the input $\bm{d}$ has no bipartite realizations.
Note that the main focus of \cite{Wang2021} is to solve a decision problem without trying to gererate any bipartite
realizations. If the goal is to generate all bipartite realizations of the input $\bm{d}$ without caring about any non-bipartite
triangle-free realizations, then the decision algorithm of \cite{Wang2021} can help our extension algorithm for bipartite
realizations skip stage 2 for many more possible inputs $\bm{d}$ compared to our generation algorithm for triangle-free realizations
described in Section \ref{sec:our_alg_tf}. We denote the set of all labeled bipartite realizations of $\bm{d}$ as $\mathcal{B}(\bm{d})$.

\section{Complexity Analysis}
\label{sec:analysis}
We just give an informal discussion of the complexity of our algorithm for triangle-free realizations in Section \ref{sec:our_alg_tf}.
Stage 1 of the algorithm is simple and it is easy to see that it takes $O(n^2)$ time. Stage 2 is complex and it is hard
to give a rigorous and accurate runtime complexity analysis for each specific case. However, it is no higher than that
of the general algorithm for all labeled realizations in \cite{Wang2024} and it is often much lower due to the applied
prunings for those graphical sequences whose majority number of realizations are not triangle-free. The best-case
runtime complexity is $O(n)$, for example on the sequence $((n-1)^1,1^{n-1})$, which is the degree sequence of
the star graph of order $n$. We anticipate it would be very
challenging to give an average-case runtime complexity across all zero-free graphical sequences of length $n$
or across the subset of all such sequences that would force the algorithm to enter stage 2. The space complexity
is easily seen to be $O(n^2)$.

We conjecture that the decision problem $P(\bm{d},3)$ might be a candidate problem in NP that is neither NP-complete
nor in P. That is, this may be an NP-intermediate problem if $P\ne NP$, as
conjectured also for the two well-known problems Graph-Isomorphism and Integer Factorization \cite{AroraBarak2009}.

\section{Experimental Evaluations}
\label{sec:experiments}
We implemented our exhaustive generation algorithm in Section \ref{sec:our_alg_tf} for triangle-free realizations
and the extension generation algorithm in \ref{sec:our_alg_bp} for bipartite realizations in C++ and
the data structure we chose for the labeled graph $G$ realizing the input $\bm{d}$ is adjacency lists.

We first use the following strategy to verify the correctness of our algorithm design and implementation.
We use an exhaustive generation algorithm from \cite{Ruskey1994} to generate all zero-free graphical
sequences of length $n$. For each such graphical sequence $\bm{d}=\{d_1\ge d_2\ge \cdots\ge d_n\}$,
we use our implementation of the exhaustive generation algorithms to generate and
count the number of labeled triangle-free (resp. bipartite) graph realizations $|\mathcal{S}(\bm{d})|$
(resp. $|\mathcal{B}(\bm{d})|$). Suppose $\bm{d}$ can also be written as $b_1^{c_1},b_2^{c_2},
\cdots,b_s^{c_s}$, where $b_1>b_2>\cdots>b_s$ are all the distinct degrees in $\bm{d}$ and there
are $c_i$ copies of $b_i$ in $\bm{d}$ for $1\le i\le s$. Then the number of labeled  triangle-free (resp. bipartite)
realizations of $\bm{d}$ without the requirement of the lexicographical ordering of the vertex labels matching the
ordering of the vertex degrees is $|\mathcal{S}(\bm{d})|\frac{n!}{\prod_{i=1}^{s}c_i!}$ (resp. $|\mathcal{B}(\bm{d})|
\frac{n!}{\prod_{i=1}^{s}c_i!}$). When we sum these values for all zero-free graphical sequences of length $n$,
we should obtain the total number of labeled triangle-free (resp. bipartite) graphs with
$n$ non-isolated vertices, which can also be independently calculated via OEIS A345218 which records
the number of labeled connected triangle-free graphs on $n$ vertices and A001832 which records the number
of labeled connected bipartite graphs on $n$ vertices (using a variant of Euler-transform \cite{Sloane1995}).
We used this approach to verify the correctness of our implementations for all $n\le 12$.
For $n\ge 13$, it appears that for some graphical sequences $\bm{d}$ the number of output labeled triangle-free
graphs $|\mathcal{S}(\bm{d})|$ is too large to finish in a reasonable amount of time.

Another strategy we can use to verify the correctness of our algorithm design and implementation is to use the
well-known tool \textbf{uniqg} from \textbf{nauty} \cite{McKay2014}. We can perform the following verification
steps: (1) generate all labeled triangle-free and bipartite  realizations of $\bm{d}$ with our implementations; (2)
obtain the number of unlabled triangle-free and bipartite realizations of $\bm{d}$ with the tool \textbf{uniqg};
and (3) sum these values over all zero-free graphical sequences $\bm{d}$ of length $n$ to get
the number of unlabeled triangle-free and bipartite graphs with $n$ non-isolated vertices, which can also be
independently calculated (using the Euler-transform \cite{Sloane1995}) via OEIS A024607 which records
the number of connected triangle-free graphs on $n$ unlabeled vertices and A005142 which records
the number of connected bipartite graphs on $n$ unlabeled vertices. We also used this approach to verify the
correctness of our implementations for all $n\le 12$.

In Table \ref{tbl:Num_TFBP_Realizations} we show some computational results on specific input instances $\bm{d}$
together with the number $|\mathcal{G}(\bm{d})|$ of all labeled realizations, the number $|\mathcal{S}(\bm{d})|$
of all labeled triangle-free realizations, and the number $|\mathcal{B}(\bm{d})|$ of all labeled bipartite realizations.
Some of these example inputs also appear in \cite{Wang2024}. Each of these results for $|\mathcal{S}(\bm{d})|$
or $|\mathcal{B}(\bm{d})|$ can be computed with our implementations in at most five seconds.

\begin{table}[!htb]
	\centering
	\caption{Some numerical results of $|\mathcal{S}(\bm{d})|$ and $|\mathcal{B}(\bm{d})|$}
	\begin{tabular}{|c|c|c|c|}
		\hline
		$\bm{d}$ & $|\mathcal{G}(\bm{d})|$ & $|\mathcal{S}(\bm{d})|$ & $|\mathcal{B}(\bm{d})|$\\
		\hline
		2,2,2,2,2 & 12 & 12 & 0\\
		\hline
		3,3,2,2,2 & 7 & 1 & 1\\
		\hline
		3,3,2,2,1,1 & 17 & 6 & 6\\
		\hline
		4,3,3,2,2,1,1 & 65 & 5 & 5\\
		\hline
		4,3,3,3,3,3,3 & 810 & 0 & 0\\ 
		\hline
		4,4,3,3,3,2,2,1 & 1931 & 33 & 33 \\
		\hline
		3,3,3,3,3,3,3,3 & 19355 & 3360 & 840\\
		\hline
		4,4,4,4,3,3,3,3 & 14634 & 36 & 36\\ 
		\hline
		6,4,4,4,4,4,3,3,2 & 109055 & 0 & 0\\ 
		\hline
		4,4,4,4,3,3,3,3,3,3 & 17262690 & 203220 & 27180\\ 
		\hline
		6,5,4,4,3,3,3,2,2,2 & 798279 & 94 & 87\\ 
		\hline
		3,3,3,3,3,3,3,3,3,3 & 11180820 & 1829520 & 257040\\ 
		\hline
	\end{tabular}
	\label{tbl:Num_TFBP_Realizations}
\end{table}
Some rare graphical sequences $\bm{d}$ satisfy $\Omega(\bm{d})=2$ so all their realizations are triangle-free.
For example, $(3^1,1^3)$ and $(2^5)$. 
The majority of graphical sequences $\bm{d}$ satisfy $\Omega(\bm{d})>2$ so not all of their realizations are triangle-free.
Some graphical sequences $\bm{d}$ satisfy $\omega(\bm{d})>2$ so none of their realizations is triangle-free.
For example, $(4^1,3^6)$ and $(6^1,4^5,3^2,2^1)$. For the latter sequence our algorithm can conclude that there is
no triangle-free realization without entering stage 2. For the former sequence our algorithm cannot conclude that there is
no triangle-free realization in stage 1 so it has to enter stage 2 to traverse the entire search tree to report this fact.

Among those graphical sequences $\bm{d}$ with at least one triangle-free realization, there are three possibilities regarding
how many of its triangle-free realizations are bipartite:
(1) all of the triangle-free realizations are bipartite (e.g. the graphical sequences $(3^2,2^2,1^2)$ and $(4^4,3^4)$);
(2) some but not all of the triangle-free realizations are bipartite (e.g. the graphical sequences $(4^4,3^6)$ and $(3^{10})$); or
(3) none of the triangle-free realizations are bipartite (e.g. the graphical sequences $(2^5)$ and $(8^3,7^1,4^7,3^1,2^1)$
with the latter sequence not shown in Table \ref{tbl:Num_TFBP_Realizations}). Given such a graphical sequence $\bm{d}$
(i.e. $\omega(\bm{d})=2$), surely we can use the extension generation algorithm in Section \ref{sec:our_alg_bp}
to tell which of the above three situations will happen. However, it may still be challenging to decide which of the
three would occur if a different approach is to be used.

Because of the effective pruning techniques employed in our algorithm, our implementation can often finish
rather quickly and much faster than that to generate all labeled realizations in \cite{Wang2024}. For example,
on the sequence $(7^1,5^9,3^1,2^1,1^1 )$ of length 13 it can finish in about 5 seconds to find all 27720 labeled
 triangle-free realizations and all 17640 labeled bipartite realizations, while it will take a long time
($>35$ days) to generate all its labeled realizations. 

We also adapted our previous C++ implementation of the James-Riha algorithm in \cite{JamesRiha1976} to generate all realizations
of a given graphical sequence $\bm{d}$ mentioned in \cite{Wang2024} by incorporating the pruning techniques
introduced in Section \ref{sec:our_alg_tf} and Section \ref{sec:our_alg_bp} and generating only triangle-free or
bipartite realizations. In order to verify that our adapted implementation of this algorithm is correct, we filter out
the isomorphic duplicates in the output graphs using the tool \textbf{uniqg} and get the number of unlabeled
triangle-free or bipartite realizations for some short graphical sequences $\bm{d}$ with $n\le 17$ whose numbers
of such unlabeled realizations have been pre-calculated.
All the computational results we obtained in our experiments match those pre-calculated values.
The number of output graphs from our adapted algorithm for an input $\bm{d}$ is often smaller than
$|\mathcal{S}(\bm{d})|$ or $|\mathcal{B}(\bm{d})|$, although it is usually still larger than the number of unlabeled
triangle-free or bipartite realizations of $\bm{d}$. This is because the adapted algorithm intends to
reduce the number of isomorphic duplicates in the outputs but cannot guarantee to remove them completely.
Also, the runtime is often shorter than that of the implementation of our exhaustive generation algorithms in
Section \ref{sec:our_alg_tf} and Section \ref{sec:our_alg_bp} because the two algorithms are
not accomplishing exactly the same task.

\section{Conclusion}
\label{sec:conclusion}
In this paper we presented an algorithm to generate all labeled triangle-free realizations of an arbitrary graphical sequence
$\bm{d}$ in lexicographical order. The generated graphs can be further tested to determine if they are bipartite to
generate all labeled bipartite realizations. The algorithm can be used to solve the decision problem $P(\bm{d},3)$ of
whether an input graphical sequence $\bm{d}$ has any triangle-free realizations. However, the complexity of deciding
$P(\bm{d},3)$ or the more general problem $P(\bm{d},r)$ remains unknown.



\begin{thebibliography}{00}
\bibitem{AroraBarak2009}
Sanjeev Arora and Boaz Barak.
\newblock Computational Complexity: A Modern Approach.
\newblock {\em Cambridge University Press}, 2009.

\bibitem{Bauer2008}
Douglas Bauer, S. Louis Hakimi, Nathan Kahl, E. Schmeichel.
\newblock Best Monotone Degree Bounds for Various Graph Parameters.
\newblock {\em Congressus Numerantium}, 192, 2008.

\bibitem{ErdosCallai1960}
P. Erdős and T. Gallai.
\newblock Graphs with given degree of vertices.
\newblock {\em Mat. Lapok}, 11:264--274, 1960.

\bibitem{Favaron1991}
O. Favaron, M. Mah{\'e}o and J.-F. Sacl{\'e}.
\newblock On the residue of a graph.
\newblock {\em Journal of Graph Theory}, 15(1):39--64, 1991.

\bibitem{Goldreich2008}
Oded Goldreich.
\newblock Computational Complexity: A Conceptual Perspective.
Cambridge University Press, 2008.

\bibitem{Hakimi1962}
S. L. Hakimi.
\newblock On Realizability of a Set of Integers as Degrees of the Vertices of a Linear Graph. {I}.
\newblock {\em Journal of the Society for Industrial and Applied Mathematics}, 10(3):496--506, 1962.

\bibitem{Havel1955}
V. Havel.
\newblock A remark on the existence of finite graphs.
\newblock {\em \v{C}asopis pro p\v{e}stov\'{a}n\'{i} matematiky}, 4:477--480, 1955.

\bibitem{JamesRiha1976}
K. R. James and W. Riha.
\newblock Algorithm 28 algorithm for generating graphs of a given partition.
\newblock {\em Computing}, 16:153--161, 1976.

\bibitem{Jelen1996}
Frank Jelen.
\newblock $k$-independence and $k$-residue of graphs with prescribed degree sequence (in German).
\newblock {\em Diploma Thesis}, University Bonn, 1996.

\bibitem{Karp1972}
Richard M. Karp.
\newblock Reducibility among combinatorial problems.
\newblock {\em Complexity of Computer Computations}, 85--103, 1972.

\bibitem{Maffray1996}
Frédéric Maffray and Myriam Preissmann.
\newblock On the NP-completeness of the $k$-colorability problem for triangle-free graphs.
\newblock {\em Discrete Mathematics}, 162(1-3):313--317, 1996.

\bibitem{Mantel1907}
W. Mantel.
\newblock Problem 28 (Solution by {H}. {Gouwentak}, {W}. {Mantel}, {J}. {Teixeira} de {Mattes}, {F}. {Schuh} and {W}. {A}. {Wythoff}).
\newblock {\em Wiskundige Opgaven}, 60--61, 1907.

\bibitem{McKay2014}
B. D. McKay and A. Piperno.
\newblock Practical Graph Isomorphism, II.
\newblock {\em Journal of Symbolic Computation}, 60:94--112, 2014.

\bibitem{Murphy1991}
O. Murphy.
\newblock Lower bounds on the stability number of graphs computed in terms of degrees.
\newblock {\em Discrete Mathematics}, 90(2):207--211, 1991.

\bibitem{Mycielski1955}
J. Mycielski.
\newblock Sur le coloriage des graphes.
\newblock {\em Colloq. Math.}, 3(2):161--162, 1955.

\bibitem{Ruskey1994}
Frank Ruskey, Robert Cohen, Peter Eades and Aaron Scott.
\newblock Alley {CATs} in search of good homes.
\newblock In {\em 25th S.E. Conference on Combinatorics, Graph Theory, and Computing}, 102:97--110. Congressus Numerantium, 1994.

\bibitem{Sierksma1991}
Gerard Sierksma and Han Hoogeveen.
\newblock Seven Criteria for Integer Sequences Being Graphic.
\newblock {\em Journal of Graph Theory}, 15(2):223--231, 1991.

\bibitem{Sloane1995}
N. J. A. Sloane and S. Plouffe.
\newblock The Encyclopedia of Integer Sequences.
San Diego, CA: Academic Press, pp. 20--21, 1995.

\bibitem{WangKleitman1973}
D. L. Wang and D. J. Kleitman.
\newblock On the {Existence} of $n$-{Connected} {Graphs} with {Prescribed} {Degrees} ($n\ge 2$).
\newblock {\em Networks}, 3(3):225--239, 1973.

\bibitem{Wang2021}
Kai Wang.
\newblock An Efficient Algorithm to Test Potential Bipartiteness of Graphical Degree Sequences.
\newblock {\em Theory and Applications of Graphs}, 8(1), 2021.

\bibitem{Wang2024}
Kai Wang.
\newblock An Efficient Algorithm to Generate all Labeled Graphs with a given Graphical Degree Sequence.
\newblock {\em Computational Science and Computational Intelligence}, vol 2506, Springer, Cham, 2024.

\bibitem{Yin2012}
J. H. Yin.
\newblock A short constructive proof of {A}.{R}. {R}ao's characterization of potentially $K_{r+1}$-graphic sequences.
\newblock {\em Discrete Applied Mathematics}, 160(3):352--354, 2012.
\end{thebibliography}
\end{document}